\documentclass[5p, number,sort&compress,times,fleqn]{elsarticle}

\usepackage{amsmath}
\usepackage{amssymb}
\usepackage{graphicx}
\usepackage{subfig}
\usepackage{bm}
\usepackage[colorlinks = true,
            linkcolor = blue,
            urlcolor  = blue,
            citecolor = blue,
            anchorcolor = blue]{hyperref}

\usepackage{tikz}
\usetikzlibrary{matrix}

\usepackage[noend]{algorithmic}
\usepackage{algorithm}

\renewcommand{\vec}[1]{\ensuremath{\bm{#1}}}
\newcommand{\ary}[1]{\boldsymbol { \mathsf{#1}} }

\DeclareMathOperator{\dif}{d}

\DeclareMathOperator*{\argmin}{arg\,min}

\newcommand{\trans}{{\mathsf{T}}}

\newcounter{ALC@tempcntr}
\newcommand{\LCOMMENT}[1]{%
 	\setcounter{ALC@tempcntr}{\arabic{ALC@rem}}
     \setcounter{ALC@rem}{1}
     \item \textit{\small{// #1}}
     \setcounter{ALC@rem}{\arabic{ALC@tempcntr}}
 }%

\usepackage[normalem]{ulem}
\newcommand{\REMOVEDTEXTMODE}[1]{\textcolor{blue}{\sout{#1}}}
\newcommand{\REMOVED}[1]{
\ifmmode
  \text{\REMOVEDTEXTMODE{$#1$}}
\else
  \REMOVEDTEXTMODE{#1}
\fi}
\newcommand{\ADDEDTEXTMODE}[1]{\textcolor{black}{#1}}
\newcommand{\ADDED}[1]{
\ifmmode
  \text{\ADDEDTEXTMODE{$#1$}}
\else
  \ADDEDTEXTMODE{#1}
\fi}

\definecolor{darkolive}{rgb}{0.01, 0.75, 0.24}
\newcommand{\SADDEDTEXTMODE}[1]{\textcolor{black}{#1}}
\newcommand{\SADDED}[1]{
\ifmmode
  \text{\SADDEDTEXTMODE{$#1$}}
\else
  \SADDEDTEXTMODE{#1}
\fi}

\begin{document}

\begin{frontmatter}

\title{Isogeometric shape optimisation of shell structures using multiresolution subdivision surfaces} 

\author{Kosala Bandara}  
\author{Fehmi Cirak\corref{cor1}}  

\cortext[cor1]{Corresponding author}

\address{Department of Engineering, University of Cambridge, Trumpington Street, Cambridge CB2 1PZ, U.K.}

\begin{abstract}
We  introduce the isogeometric shape optimisation of thin shell structures using subdivision surfaces. Both triangular Loop  and quadrilateral Catmull-Clark subdivision schemes are considered for geometry modelling and  finite element analysis. A gradient-based  shape optimisation technique is implemented to minimise compliance, i.e.  to maximise stiffness.   Different control meshes describing  the same surface are used for geometry representation, optimisation and finite element analysis.  The finite element analysis is performed with subdivision basis functions corresponding to a sufficiently refined control mesh. During iterative shape optimisation the geometry is updated starting from the coarsest control mesh and  proceeding to increasingly finer control meshes. This  multiresolution approach provides a  means for regularising the optimisation problem and prevents the appearance of sub-optimal jagged geometries with fine-scale oscillations.  The finest  control mesh for optimisation is chosen in accordance with the desired smallest feature size in the optimised geometry.   
The  proposed approach is applied to three optimisation examples, namely a catenary, a roof over a rectangular domain and a freeform architectural shell roof. The influence of the geometry description and the used subdivision scheme on the  obtained optimised curved geometries is investigated in detail. 
\end{abstract}

\newpage


\end{frontmatter}

\section{Introduction}
%

Shell structures are curved solids with one dimension significantly smaller than the other two. They are prevalent in many engineering applications, most prominently in aerospace, automotive and structural engineering. The load carrying capacity of shells can be greatly increased by systematically optimising their curved shape.  Due to their small thickness the mechanics of shells can be efficiently described with surface models. The mechanical response of a thin shell depends, according to the Kirchhoff-Love model, on the  first and second fundamental forms of the surface. In shape optimisation of shells,  the efficient and flexible description of  freeform surfaces and the finite element discretisation of the governing equations defined on them are intrinsically linked.  In this paper we use the subdivision surfaces  as a common representation for geometric modelling and finite element discretisation of  Kirchhoff-Love shell equations.

\begin{figure*}
	\centering 
  	\begin{minipage}[b]{0.33\textwidth} 
  		\subfloat[][Coarse control mesh.]{\includegraphics[width=\textwidth]{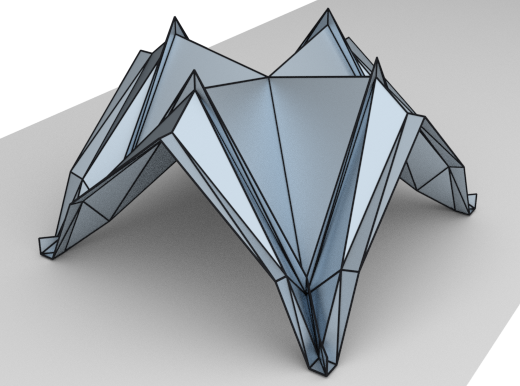}}
	\end{minipage}
	  	\begin{minipage}[b]{0.33\textwidth} 
  		\subfloat[][Finite element mesh.]{\includegraphics[width=\textwidth]{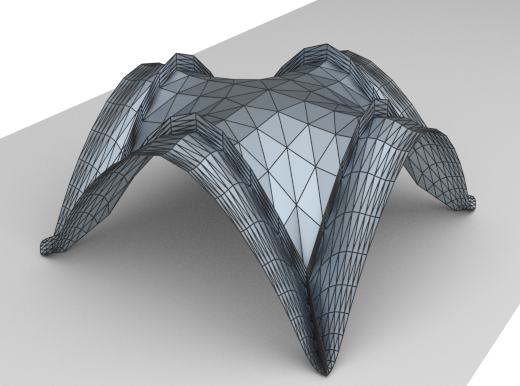}}
	\end{minipage}
  	\begin{minipage}[b]{0.33\textwidth} 
  		\subfloat[][Vertical displacement isocontours. ]{\includegraphics[width=\textwidth]{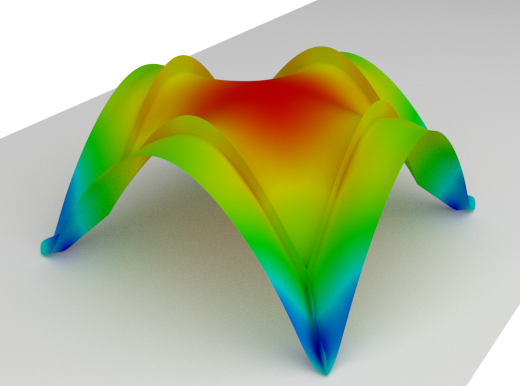}}
	\end{minipage}
  	\caption{Shape optimised thin shell roof structure emanating from a flat plate with stiffeners. The structure is loaded with a uniformly distributed vertical load and is supported at its four corners.  The vertex positions of the coarse control mesh (a) represent the optimisation design variables. The compliance cost function and its derivatives are evaluated with  finite element analysis using the twice subdivided control mesh (b) and the corresponding smooth subdivision basis functions. The limit surface and the isocontours of the vertical displacement are shown in (c).  }
  \label{fig:introExample}
\end{figure*}

Isogeometric analysis aims to unify geometric  modelling and finite element analysis by using for the latter usual computer-aided design (CAD)  basis functions, like NURBS. Since its inception by Hughes et al.~\cite{Hughes:2005aa,Cottrell:2009aa} isogeometric analysis has become immensely popular and has been applied to a wide range of engineering problems, too many to list here. Prior to the advent of isogeometric analysis, the integrated geometric modelling and finite element analysis of shells using subdivision surfaces was proposed in~\cite{Cirak:2000aa}.  Specifically, Loop subdivision surfaces were used for discretising the Kirchhoff-Love shells  and representing their geometry. As an extension of this approach, the treatment of  industrially prevalent non-manifold shell geometries and the inclusion of out-of-plane shear deformations relevant for thicker shells were proposed in~\cite{Cirak:2011aa}  and~\cite{Long:2012aa}, respectively. More recently the isogeometric analysis of shells and beams using NURBS basis functions were introduced in~\cite{Kiendl:2009aa,nagy2010isogeometric}. The use of  smooth subdivision and NURBS   basis functions has also the advantage that they have square-integrable curvatures, which is necessary for discretising the Kirchhoff-Love shell equations depending on curvatures.

 In the present work, we investigate  the gradient-based shape optimisation of  shell structures using  subdivision surfaces for geometric modelling and finite element analysis. The minimised cost function is the compliance so that (qualitatively)   displacements, strains and stresses are minimised.  In a typical industrial design setting both the input to and output from structural optimisation is a geometry, i.e.  a CAD model. During optimisation the cost function and its derivatives with respect to some geometric design parameters, i.e. design sensitivities, need to be computed with  finite element analysis~\cite{Haug:1986aa,Christensen:2009aa}. Hence, as a matter of fact, the  interoperability of geometry and finite element models is crucial. Equally important are techniques  for choosing suitable  geometric design variables that  can parameterise a sufficiently large set of geometries. Over the years,  a wide variety of shape parameterisation techniques have been proposed that are  based either on  a CAD, a finite element analysis (FEA) or an intermediary model~\cite{Haftka:1986aa, samareh2001survey}.  In  one group of  techniques  the geometric design variables are the parameters of the CAD model or a reconstructed CAD-like spline model~\cite{Braibant:1984aa,bletzinger1993form,cervera2005evolutionary,robinson2012optimizing,han2014adaptive}. In the second group of techniques the design variables are  the vertex positions of the finite element mesh~\cite{le2011gradient,firl2013regularization,bletzinger2014consistent}.  Yet in another group of parameterisation  techniques, like the ones based on  radial basis functions~\cite{jakobsson2007mesh} or free-form deformations~\cite{imam1982three, sederberg1986free}, the design variables are only indirectly linked to the CAD or FEA model. This list of  parameterisation techniques is not intended to be complete. The abundance of shape parameterisation techniques is partly due to the inherent limitations and incompatibilities of conventional CAD and FEA models in the context of shape optimisation. With isogeometric analysis  using subdivision surfaces most of the incompatibilities between the CAD and FEA  representations can be elegantly circumvented. To this end, the increasing availability of subdivision surfaces in CAD systems, including PTC Creo, CATIA, Siemens NX or  Autodesk Fusion 360,  is noteworthy.

In isogeometric shape optimisation with subdivision surfaces different resolutions, i.e. control meshes, of a surface are employed for optimisation and  performing the finite element analysis, see Figure~\ref{fig:introExample}. In computer graphics subdivision surfaces are usually viewed as a process for generating increasingly finer meshes that converge in the limit to a surface~\cite{Zorin:2000aa}. Alternatively, they can be viewed as a generalisation of splines to arbitrary connectivity meshes~\cite{Peters:2008aa}. In the proposed  optimisation approach, both viewpoints are simultaneously exploited. Subdivision surfaces are best considered as generalised splines when used as finite element basis functions. On the other hand, the discrete computer graphics viewpoint with the associated data structures and algorithms is best suited  for simultaneously operating on different resolutions in a memory and time efficient manner. 
Our present implementation is based on the triangular Loop~~\cite{Loop:1987aa} and quadrilateral Catmull-Clark schemes~\cite{Catmull:1978aa}, or more specifically on their extended versions introduced in~\cite{Biermann:2000aa}. The finite element analysis is performed with basis functions corresponding to a sufficiently fine control mesh. Within the  optimisation loop, starting with the coarsest,  the vertex positions of increasingly finer control meshes  are used as design variables. The resolution of the control mesh determines the extent of applied geometry changes, because each vertex has control over  the surface within a two-ring of adjacent elements.  The derivatives of the cost function with respect to vertex positions is first computed on the fine finite element control mesh and subsequently projected to the coarser control meshes corresponding to the design variables. This projection provides a means for smoothing, or filtering, of the computed design sensitivities and  prevents the appearance of jagged optimised geometries with fine-scale oscillations. The need for such a smoothing, or filtering,  in shape optimisation is widely discussed in  literature, see e.g.~\cite{le2011gradient,bletzinger2014consistent} and references therein.

An earlier two-level version of the multiresolution optimisation approach proposed in  this paper was introduced in~\cite{Cirak:2002aa}. In that exploratory work Loop subdivision and non-gradient based  optimisation algorithms were used.  More recently the  proposed multiresolution approach has been applied to other types of optimisation problems, namely electrostatic  shape optimisation  of high-voltage devices~\cite{Bandara:2014aa} and  shape optimisation of volumetric solids~\cite{bandara2016shape}. The electrostatic simulations are performed with the boundary element method and the solid simulations with the voxel-based immersed finite element method.

This paper is organised as follows. In Section~\ref{sec:thinShells} we begin by reviewing the Kirchhoff-Love  model for thin shells and its discretisation with subdivision basis functions.  We then introduce the compliance optimisation problem and compute with an adjoint approach the cost function  derivatives with respect to the vertex positions of the subdivision control mesh. Section~\ref{sec:subdivision} provides a very brief introduction to subdivision surfaces. Subsequently, in Section~\ref{sec:optimisation} the multiresolution algorithm is introduced. Finally, in Section~\ref{sec:examples} we introduce three examples of increasing complexity, namely the shape optimisation of a thin-strip, a roof over a rectangular domain and a freeform architectural shell roof. 

\section{Thin-shells \label{sec:thinShells}}
%
\subsection{Governing equations and discretisation \label{sec:governing}}
\begin{figure*}
	\centering
  \includegraphics[scale=0.9]{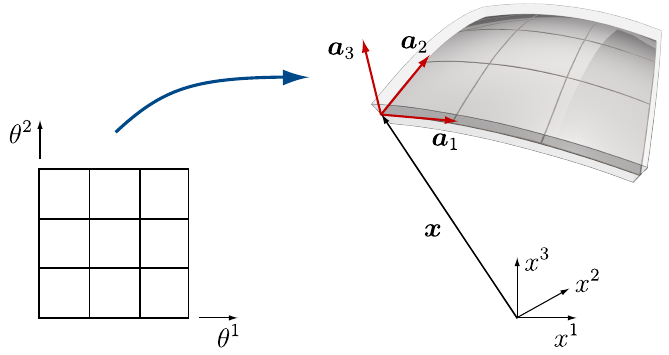}
  \caption{Solid thin-shell and its mid-surface (right) with the position vector~$\vec x(\theta^1, \theta^2)$, the covariant basis vectors~$\vec a_1$ and~$\vec a_2$, and the normal vector~$\vec a_3$. The parameter space with the coordinates~$(\theta^1, \theta^2)$ is shown on the left.}
  \label{fig:thinShell}
\end{figure*}

%
We consider a thin-shell with the undeformed  mid-surface~$\Omega$ and the thickness~$t$, see Figure \ref{fig:thinShell}. The surface~$\Omega$ is parameterised with the curvilinear coordinates~$(\theta^1, \theta^2) \in \varmathbb{R}^2$  providing each material point on the surface with  a unique parametric coordinate. The position of the material points is denoted with the coordinates \mbox{$\vec x(\theta^1, \theta^2) \in \varmathbb{R}^3$}.

The standard covariant basis vectors of the mid-surface~$\vec a_\alpha$ and the unit normal~$\vec a_3$ are given by 
\begin{equation}
	\vec a_\alpha = \frac{ \partial \vec x}{\partial \theta^\alpha } = \vec x_{,\alpha} \, , \quad  \vec a_{3} = \frac{\vec a_1 \times \vec a_2}{| \vec a_1 \times \vec a_2 |}  \, .
\end{equation}
The corresponding contravariant basis vectors~$\vec a^\alpha$  are defined through the relation $\vec a^\alpha \cdot \vec a_\beta = \delta^\alpha_\beta$, where~$ \delta^\alpha_\beta$ is the Kronecker delta. Here and in the following  the Greek indices take the values~$\{ 1,2\} $ and the summation convention is used.  
 
Displacing each material point on the mid-surface with a displacement vector $\vec u (\theta^1, \theta^2) \in \varmathbb{R}^3$ yields a deformed (or, displaced) mid-surface. Subject to few mechanical assumptions, it can be shown that the differences in the first and second fundamental forms of the original and displaced surface provide  suitable strain measures.  The difference in the first fundamental form is referred to as the membrane strain tensor~$\vec \alpha$ and the difference in the second fundamental form as the bending strain tensor~$\vec \beta$.  In case of small displacements, as derived e.g. in~\cite{Cirak:2000aa}, the linearised membrane strain tensor is 
\begin{equation} \label{memstrain_lin}
	\vec \alpha  = \frac{1}{2} \left ( {\vec{a}}_\alpha \cdot \vec{u}_{, \beta} + \vec{u}_{, \alpha} \cdot {\vec{a}}_{\beta} \right ) \vec a^\alpha \otimes \vec a^\beta 
\end{equation}
and the linearised bending strain tensor is 
\begin{equation} \label{benstrain_lin}
	\begin{aligned} 
		\vec \beta  &= \left (- \vec{u}_{, \alpha \beta} \cdot  {\vec{a}}_3 + \frac{1}{\sqrt{ {a}}} \, [ \vec{u},_{1} \cdot(  {\vec{a}}_{\alpha , \beta} \times  {\vec{a}}_2) + \vec{u},_{2} \cdot ( {\vec{a}}_1 \times {\vec{a}}_{\alpha, \beta} ) ]  \right.
\\
	 & \left. +  \frac{ {\vec{a}}_3  \cdot  {\vec{a}}_{\alpha, \beta}}{ \sqrt{a}}  \, [ \vec{u},_{1} \cdot ( {\vec{a}}_2  \times  {\vec{a}}_3) + \vec{u},_{2} \cdot ( {\vec{a}}_3 \times  {\vec{a}}_1) ]  \right )  \vec a^\alpha \otimes \vec a^\beta \,  
\end{aligned}
\end{equation}
with  $\sqrt{a} =  | \vec a_1 \times  \vec a_2| $.

Next, we consider the potential energy of the displaced shell 
\begin{equation}	\label{eq:energyKL}
\begin{aligned}
	  \Pi(\vec u)  & = \int_{\Omega}  \left ( W^m( {\vec \alpha}  ) + W^b( {\vec \beta} )  \right )  \dif \Omega 
	  				 - \int_\Omega \vec p \cdot \vec u  \dif \Omega \\ & - \int_\Gamma  \vec r \cdot \vec u \dif \Gamma ,
\end{aligned}	  
\end{equation}
where the first integral is the internal potential energy consisting of the sum of the internal membrane and bending energy densities~$W^m$ and~$W^b$, respectively. The remaining two integrals represent the external potential energy resulting from the prescribed surface load vector~$\vec p$ and the edge load vector~$\vec r$.  For an elastic material the two internal energy densities are defined with 
\begin{equation}
\begin{aligned}
	W^m (  \vec \alpha ) &= \frac{1}{2} \frac{E t }{1 - \nu^2}  \vec \alpha : \vec H : \vec \alpha  \qquad 
	\\
	W^b( \vec \beta ) &= \frac{1}{2} \frac{E t^3}{12 (1 - \nu^2)} \vec \beta : \vec H : \vec \beta \, , 
\end{aligned}	
\end{equation}
where~$E$ is  the Young's modulus, $\nu$ is the Poisson's ratio and $\vec H = H^{\alpha \beta \gamma \delta} \vec a_\alpha \otimes  \vec a_\beta  \otimes \vec a_\gamma \otimes \vec a_\delta $ is  an auxiliary fourth order tensor  with the contravariant components 
\begin{equation}
		H^{\alpha \beta \gamma \delta}  =   \nu \, {a}^{\alpha\beta} {a}^{\gamma \delta}  + \frac{1}{2} ( 1 -\nu ) \, ( {a}^{\alpha \gamma} {a}^{\beta \delta} + {a}^{\alpha \delta} {a}^{\beta \gamma} ) 
\end{equation}
and the contravariant metric~$a^{\alpha \beta} = \vec a^\alpha \cdot \vec a^\beta$.

The equilibrium configurations of the shell with prescribed loading are obtained from minimising~\eqref{eq:energyKL}. Note that for a well-posed problem also the displacements on some parts of the boundary have to be prescribed in addition to the loading. In a finite element approximation, the mid-surface position and the displacement vectors in~\eqref{eq:energyKL} are approximated with basis functions and their coefficients
\begin{equation} \label{eq:interpolation}
	\vec{x}  (\theta^1, \theta^2)  \approx  \sum_{i} N_i(\theta^1, \theta^2)\vec{x}_i  \, , \quad  	\vec{u}  (\theta^1, \theta^2)  \approx \sum_{i} N_i(\theta^1, \theta^2)\vec{u}_i  \, . 
\end{equation}
In our implementation, the basis functions~$N_i(\theta^1, \theta^2) $ are obtained either from  triangular Loop subdivision or quadrilateral Catmull-Clark subdivision. In both schemes there is one basis function associated with each vertex of the control mesh. Hence, the coefficients~$\vec x_i$  and~$\vec u_i$ are simply the position and displacements of a (control) vertex with the index~$i$. Introducing the approximations~\eqref{eq:interpolation}  into the potential~\eqref{eq:energyKL} yields a discrete minimisation problem for computing the vertex displacements~$\vec u_i$, 
\begin{equation} \label{eq:minimisation}
	\vec u_i = \argmin_{\vec u_i} \Pi (\vec u_i)  \;  \;   \Rightarrow  \; \; \frac{\partial \Pi (\vec u_i)}{\partial \vec u_i} = \vec 0 \, .
\end{equation}
In order to compute the stationary points of~$\Pi(\vec u_i)$  domain integrals are numerically evaluated in a usual finite element fashion by iterating over the elements/faces in the control mesh. 
Around extraordinary vertices subdivision surfaces consist of  an infinite sequence of  ever smaller rings of box splines in Loop subdivision and b-splines in Catmull-Clark subdivision~\cite{Peters:2008aa}.  Hence, their numerical integration requires special care and has been investigated in several recent numerical studies~\cite{wawrzinek2016integration, juttler2016numerical}.  For practical computations,  the integration of each finite element using  Gauss integration with~$3$ points for Loop subdivision and~$4$ points for Catmull-Clark subdivision  appears to  provide the best trade-off between accuracy and robustness~\cite{Long:2012aa, juttler2016numerical}.  As an aside, the issue of accuracy of quadrature is independent from the sub-optimal convergence of finite elements based on subdivision surfaces, which is presently a very active area of research, see~\cite{majeedCirak:2016} and the references therein. At the Gauss points, we evaluate the  basis functions with a simplified version of the algorithm proposed by Stam~\cite{Stam:1998aa,Stam:1999aa}, see~\cite{Cirak:2000aa,Cirak:2011aa}. Specifically, since the Gauss points are relatively far from extraordinary vertices there are no efficiency gains from the eigendecomposition considered in~\cite{Stam:1998aa,Stam:1999aa}.  After numerical integration the stationarity condition for the minimisation problem~\eqref{eq:minimisation} yields a discrete system of equations 
\begin{equation} \label{eq:discreteEquilibrium}
	   \frac{\partial \Pi (\vec u_i)}{\partial \vec u_i} = \vec 0     \;  \;   \Rightarrow  \; \; \ary K \ary u = \ary f, 
\end{equation}
where~$\ary K$ is the symmetric, positive-definite system (or, stiffness) matrix,~$\ary u$ is the  array of  vertex displacements containing all~$\vec u_i$ and~$\ary f$ is the array of corresponding vertex forces. For further details see~\cite{Cirak:2000aa,Cirak:2002aa}.

\subsection{Design sensitivities}
%
In shape optimisation we consider a shell structure with prescribed loading and displacement boundary conditions and  aim to find its mid-surface such that a user chosen cost function  
\begin{equation}
	\min_{\ary x } J (\ary x, \ary u) 
\end{equation}
is minimised. As a constraint the array of the vertex displacements~$\ary u$ has to satisfy the equilibrium equations~\eqref{eq:discreteEquilibrium}. In practice, there are additional constraints, for instance pertaining to the surface area of the shell or the position of selected vertices, which will be omitted in this section. Moreover, in all examples presented in this paper the cost function is the compliance of the structure 
\begin{equation} \label{eq:discreteCost}
	J (\ary x, \ary u) =  \ary u^\trans \ary f =  \ary u^\trans \ary K \ary u \, .
\end{equation}
Informally, minimising the compliance leads to stiffer shell structures with smaller displacements~$\ary u$. The subsequent derivations carry over to other cost functions, see e.g.~\cite{Haug:1986aa}.  

In order to use a gradient-based optimisation algorithm for minimising~$J (\ary x, \ary u) $ the  derivatives of the cost function with respect to the vertex coordinates, also referred to  as design sensitivities or shape gradients, are needed. To this end, we consider the adjoint formulation with  
\begin{equation}
  \label{eq:lagrangianDiscrete}
  L (\ary{x}, \ary{u},  \ary{\lambda}) =  J (\ary x, \ary u) +  \ary{\lambda}^\trans [ \ary{f}- \ary{K}  \ary{u}]  =  \ary{u}^\trans \ary{K} \ary u+  \ary{\lambda}^\trans [ \ary{f}- \ary{K}  \ary{u}]  ,
\end{equation}
where~$\ary \lambda$ is an array of Lagrange parameters. The stationarity condition for~$L (\ary{x}, \ary{u},  \ary{\lambda}) $ with respect to the vertex  displacements  leads to the adjoint problem 
\begin{equation}
    \frac{\partial L(\ary x, \ary u, \ary \lambda) }{\partial \ary u}  = \ary{0} \quad \Rightarrow \quad    \ary K \ary \lambda  = 2 \ary K \ary u    \quad \Rightarrow \quad  \ary \lambda = 2 \ary  u \, . \label{eq:adjointDiscrete}
\end{equation}
Here, we made use of the symmetry of the stiffness matrix~$\ary K$.  The equality between the  Lagrange parameters~$\ary \lambda$ and displacements~$\ary u$ (up to the constant~$2$) is only valid when the  cost function is the compliance~\eqref{eq:discreteCost}.  The stationarity condition for $L (\ary{x}, \ary{u},  \ary{\lambda}) $ with respect to the vertex coordinates leads to the design sensitivities 
\begin{align}
	\frac{ \partial L (\ary{x}, \ary{u},  \ary{\lambda}) }{\partial \ary x} &=  \ary{u}^\trans \frac{\partial \ary K}{\partial \ary x} \ary u + 2 \ary u^\trans \left [  \frac{\partial \ary f}{\partial \ary x} - \frac{\partial \ary K}{\partial \ary x} \ary u \right ]  \\
	&=   2 \ary u^\trans \frac{\partial \ary f}{\partial \ary x} -\ary u^\trans \frac{\partial \ary K}{\partial \ary x} \ary u \, . 
\end{align}
At equilibrium, that is when~$ \ary K \ary u - \ary f = \ary 0$ is exactly satisfied, the gradients of the Lagrangian~$L(\ary x, \ary u, \ary \lambda)$ and the cost function~$ J(\ary x, \ary u)$ with respect to the vertex coordinates are identical. Hence, in gradient-based shape  optimisation the vertex coordinates have to be  perturbed in the direction   
\begin{equation} \label{eq:finalGradient}
	\ary v =  - 2 \ary u^\trans \frac{\partial \ary f}{\partial \ary x} + \ary u^\trans \frac{\partial \ary K}{\partial \ary x} \ary u \, .  
\end{equation}
to achieve a decrease in the cost function. In order to compute~$\ary v$ the derivatives of the load vector $\ary f$ and system (or  stiffness) matrix $\ary K$ with respect to the vertex coordinates $\ary x$ are needed. For the considered Kirchhoff-Love shell formulation it is straightforward to compute these derivatives by systematic element-by-element differentiation of the discrete equilibrium equations~\eqref{eq:discreteEquilibrium}. Note that other shell formulations, especially  ones available in commercial software, contain as degrees of freedom in addition to vertex positions also the rotations of  vertex director vectors. This usually makes the computation of the related design sensitivities more complex.

\section{Subdivision surfaces \label{sec:subdivision}}
%

In the isogeometric analysis context it is expedient to consider subdivision surfaces as the generalisation of splines to arbitrary connectivity meshes~\cite{Zorin:2000aa,Peters:2008aa}.  As known, refinable basis functions allow to represent the same  spline surface  with  control meshes of different resolutions. Loop subdivision~\cite{Loop:1987aa} is the generalisation of quartic box-splines to arbitrary connectivity triangular meshes and Catmull-Clark~\cite{Catmull:1978aa}  subdivision is the generalisation of tensor-product cubic b-splines to arbitrary  connectivity quadrilateral meshes. Both schemes lead to basis functions that are refinable. 

The control meshes are refined by quadrisection of elements. In triangular Loop subdivision this is achieved by introducing a  new vertex  on each edge. In  Catmull-Clark in addition to  new vertices on the edges a new vertex at the centre of each element is created. The control vertex coordinates of a refined mesh at level~$\ell+1$ are obtained from the vertex coordinates of the coarse mesh at level~$\ell$ according to 
\begin{equation}
	\ary x^{\ell+1} = \ary{S} \ary{x}^\ell \,. 
\end{equation}
For  vertices located in the regular regions of a mesh  the subdivision matrix~$\ary S$ contains the standard knot insertion weights.  The matrix components  associated to the vertices in the irregular regions and at the boundaries are given by the specific subdivision scheme. Explicit expressions for the subdivision matrix $\ary S$ can be found, e.g., in~\cite{Zorin:2000aa}.  The successive subdivision refinement of a control mesh can be interpreted as the chain of linear mappings
\begin{equation} \label{eq:subdivRef}
	\begin{tikzpicture}[descr/.style={fill=white,inner sep=1.5pt}]
		\matrix (m) [matrix of math nodes, row sep=2.7em, column sep=2.3em, text height=1.5ex, text depth=0.25ex]	
		{ \ary{x}^0 & \ary{x}^1 &  \ary{x}^2  &   \cdots &  \ary{x}^{\ell-1}  &   \ary{x}^{\ell}  \\ };
		\path[->, font=\small]
		(m-1-1) edge node[above] {$ \ary S $} (m-1-2)
		(m-1-2) edge node[above] {$ \ary S $} (m-1-3)
		(m-1-3) edge node[above] {$ \ary S $} (m-1-4)
		(m-1-4) edge node[above] {$ \ary S $} (m-1-5)
		(m-1-5) edge node[above] {$ \ary S $} (m-1-6);
	\end{tikzpicture} .
\end{equation}
The size of the vertex coordinates~$\ary x^\ell$ increases with increasing~$\ell$ and the size of the subdivision matrix~$\ary S$ increases accordingly. All control meshes converge irrespective of their level~$\ell$  to the same surface. As mentioned in Section~\ref{sec:governing}, the algorithm proposed by Stam~\cite{Stam:1998aa,Stam:1999aa} provides a spline based parameterisation of subdivision surfaces so that the properties of arbitrary surface points can be evaluated,~cf.~\eqref{eq:interpolation}. There are also alternative parameterisations available, see~\cite{boier2004differentiable,antonelli2013subdivision,wawrzinek2016integration}.

In shape optimisation the coarsening of the refined subdivision control meshes is also needed, that is,  
\begin{equation} \label{eq:rMatrix}
	\ary x^{\ell} = \ary{R} \ary{x}^{\ell+1}  \quad \text{with }  \; \; \ary R = (\ary S^\trans \ary S)^{-1}\ary S^\trans \, , 
\end{equation}
where  the coarsening matrix~$\ary R$ is defined as  the pseudo-inverse of the subdivision matrix~$\ary S$. The span of geometries that can be represented with~$\ary x^{\ell+1}$ is larger than the ones with~$\ary x^{\ell}$, hence the specific form of~$\ary R$ is a choice. As discussed in~\cite{Bandara:2014aa, bandara2016shape},  $\ary R$ can be interpreted as a smoothing operator and accordingly different choices are possible. Similar to  subdivision refinement the coarsening matrix can be successively applied in order to obtain coarser representations of the geometry, i.e., 
\begin{equation}
	\begin{tikzpicture}[descr/.style={fill=white,inner sep=2.5pt}]
		\matrix (m) [matrix of math nodes, row sep=2.7em, column sep=2.3em, text height=1.5ex, text depth=0.25ex]	
		{ \ary{x}^0 & \ary{x}^1 &  \ary{x}^2  &   \cdots &  \ary{x}^{\ell-1}  &   \ary{x}^{\ell}  \\ };
		\path[<-, font=\small]
		(m-1-1) edge node[above] {$ \ary R $} (m-1-2)
		(m-1-2) edge node[above] {$ \ary R $} (m-1-3)	
		(m-1-3) edge node[above] {$ \ary R $} (m-1-4)
		(m-1-4) edge node[above] {$ \ary R $} (m-1-5)
		(m-1-5) edge node[above] {$ \ary R $} (m-1-6);
	\end{tikzpicture}  .
\end{equation}
Note that,  by definition~\eqref{eq:rMatrix}, one step of subdivision refinement followed by one step of coarsening does not change the vertex coordinates, or~$\ary{R} \ary{S} = \ary{I}$.

\section{Shape optimisation \label{sec:optimisation}}
%
The subdivision surfaces enable us to use different resolutions of the same geometry for optimisation and analysis.  Crucially,  in the spirit of isogeometric analysis the control meshes for analysis and optimisation represent the same surface. A simplified two-level  version of the proposed iterative gradient-based optimisation algorithm is shown Algorithm~\ref{algo:multiresTopDown}. For a more advanced multiresolution version employing a wavelet-like decomposition of the surface, in the context of shape optimisation of solids, see~\cite{bandara2016shape}. The optimisation and analysis meshes correspond to different  refinement levels in a subdivision hierarchy. The optimisation level is initialised with~$\ell_o=0$ and the analysis level is fixed with~$\ell_c=n$, where $n$ is user prescribed. The gradient of the compliance cost function~$\ary v^{\ell_c}$ is computed with a finite element analysis using basis functions of the control mesh at level~$\ell_c$, see Section~\ref{sec:thinShells}. Subsequently, through successive multiplication of the gradient  with the coarsening matrix~$\ary R$  the optimisation  level gradient~$\ary v^{\ell_o}$ is obtained.  As indicated in Algorithm~\ref{algo:multiresTopDown} the updated coordinates~$\ary x^{\ell_o}$ of the optimisation level can be, for instance, obtained with a simple steepest descent method. 
 Instead of this simple update algorithm, we use in the presented examples the MMA optimisation algorithm in the NLopt library~\cite{nloptPackage}. The MMA algorithm usually requires fewer iterations to converge and is able to consider both equality and inequality side constraints.  Moreover, in the presented examples the optimisation level is successively increased every time a stationary point is reached until a user prescribed maximum optimisation level is reached.  It is evident that~$ \ell_o \le  \ell_c$. The termination criterion for optimisation iterations is tightened with increasing optimisation level~$\ell_o$. In the presented examples the tolerance parameter for the iterations is chosen to be $\varepsilon = 10^{-3 (\ell_o+1)}$. 
\begin{algorithm}[H]
\caption{Multiresolution shape optimisation}
\label{algo:multiresTopDown}
\begin{algorithmic}[1]
\LCOMMENT{read maximum optimisation level~$\ell_{o,max}$   }
\LCOMMENT{read finite element analysis level~$\ell_c$ }
\LCOMMENT{read input coarse control mesh~$\ary x^{0}$   }

\LCOMMENT{initialise optimisation level}
\STATE $\ell_o =  0$  
\LCOMMENT{initialise cost function}
\STATE $J =  \infty$  
\LCOMMENT{iterate over optimisation levels}
\WHILE {$\ell_o \le \ell_{o,max}$}
\LCOMMENT{update vertex coordinates $\ary{x}^{\ell_o}$ while cost function decreases}
\REPEAT 
\LCOMMENT{subdivide optimisation level $\ell_o$ up to analysis level $\ell_c$}
 \FOR{$ \ell \gets \ell_o \textrm{ to } \ell_c-1$}
 	\STATE $\ary{x}^{\ell+1} \gets \ary{S} \ary x^{\ell}  $
 \ENDFOR 
\LCOMMENT{store previous cost function}
 \STATE $J_{\text{previous}}  \gets J $ 

\LCOMMENT{compute cost function $J = J (\ary x^{\ell_c}, \ary u^{\ell_c}) $ and its gradient $\ary v^{\ell_c}$}

\LCOMMENT{project gradient to optimisation level}

\FOR{$ \ell \gets \ell_c \textrm{ to } \ell_o+1$}
 	\STATE $\ary{v}^{\ell-1} \gets \ary R \ary v^{\ell}  $
 \ENDFOR 
 
\LCOMMENT{update vertex coordinates of the optimisation level}
\STATE $\ary{x}^{\ell_o}  \leftarrow  ( \ary{x}^{\ell_o} + \alpha \ary{v}^{\ell_o}) \quad
\text{ with } \alpha \ge 0$
\UNTIL{{$ ( J_{\text{previous}} -J ) <  \dfrac{ \varepsilon}{2} (J_{\text{previous}} +J)$}}

\LCOMMENT{increment optimisation level}
\STATE $\ell_o \gets (\ell_o+1)$
\ENDWHILE
\end{algorithmic}
\end{algorithm}

\section{Examples \label{sec:examples}}
%
Three examples are presented to demonstrate the functioning of the proposed isogeometric shape optimisation of thin-shell structures using subdivision. In all examples the objective is to minimise the compliance. As subdivision schemes the triangular Loop and the quadrilateral Catmull-Clark scheme are used. In order to preserve corners and edges of the original geometry, at vertices and edges on the boundary modified subdivision stencils are applied~\cite{Biermann:2000aa,Cirak:2011aa}. 

\subsection{Thin strip}
\begin{figure}
\centering
  \subfloat[][Initial geometry.]
  {\label{fig:stripInitial}
  \includegraphics[scale=0.8]{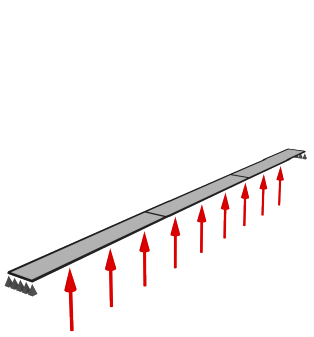}} \\
  \subfloat[][Optimised geometry.]
  {\label{fig:stripOptimised}
    \hspace{0.15\linewidth}
  \includegraphics[scale=0.8]{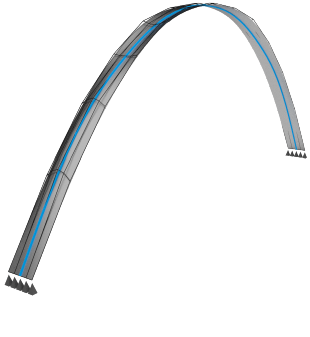}}
    \caption[]{Optimisation of an uniformly loaded thin strip.}
  \label{fig:catenaryCurveEx}
\end{figure}

%
In this verification example, a thin strip pinned at both ends and subjected to a vertical distributed load is optimised, see Figure~\ref{fig:catenaryCurveEx}. Initially, the strip is a  narrow flat plate with length~$1$ equal to the distance between the supports. The magnitude of the vertical uniformly distributed load is~$1000$. The width of the strip is~$0.05$; the thickness is~$0.02$; the Young's modulus is~$E = 2\times 10^8$; and the Poisson's ratio is~$\nu = 0.3$.  

The Catmull-Clark subdivision scheme is used for representing the geometry and finite element analysis.  Although there are no extraordinary vertices in the control mesh, the subdivision basis functions are neither uniform nor non-uniform b-splines due to the treatment of the boundaries~\cite{Biermann:2000aa}. The initial coarse mesh used for optimisation contains only~$3$ elements along the length and~$1$ element across the width of the strip. This increases to~$48$ in the twice subdivided analysis fine mesh, i.e.~$\ell_c \equiv 2$. During  compliance optimisation the mesh resolution is increased starting from~$\ell_o=0$ up to~$\ell_o =  2$.  Only the out-of-plane position of the control points in the direction of the load vector are optimised. The length of the optimised strip is chosen to be either~$1.1$, $1.2$ or $1.3$ by prescribing its area, see~\cite{bandara2016shape}  for the treatment of area constraints.

\begin{figure}
\centering
	\subfloat[][Prescribed length $=1.1$.] {
    \includegraphics[scale=0.85]{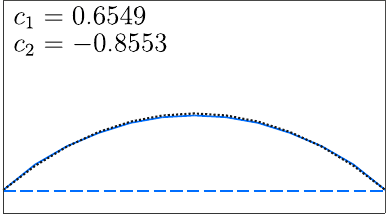}} \\
    	\subfloat[][Prescribed length $=1.2$.] {
    \includegraphics[scale=0.85]{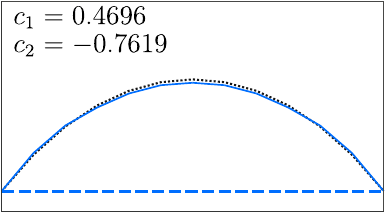}} \\
	\subfloat[][Prescribed length $=1.3$.] {
    \includegraphics[scale=0.85]{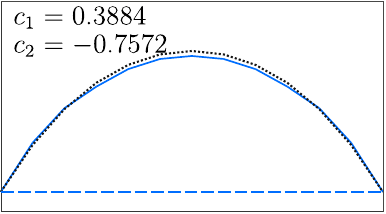}} \\
    \caption[]{Optimisation of a thin strip with different prescribed lengths. The  blue lines show the centre line of the strip before and after optimisation and the dotted black line is the catenary curve. }
  \label{fig:catenaryCurveComparison}
\end{figure}

As known from classical mechanics, the shape of the curve assumed by a loose string pinned at both ends is a catenary curve~\cite{lockwood1961book}, which is for the considered geometry of the form  
\begin{equation}
y= c_1 \cdot \cosh \left( x / c_1 \right ) + c_2 \,  , 
\label{eq:catenary}
\end{equation}
where the~$y$-axis is parallel to the applied load vector and the left and right supports have the coordinates~$(-0.5, 0)$ and~$(0.5, 0)$, respectively. The constants~$c_1$ and~$c_2$ depend on the chosen length of the optimised strip. The comparison of the optimisation results with the corresponding catenary curves is shown in Figure~\ref{fig:catenaryCurveComparison}. The reduction of the compliance cost function is more than~$99.9\%$ and the optimisation results show good visual agreement with the catenary curve. The slight deviation from the catenary is possibly due to the finite width of the strip, which leads during optimisation to some curvature generation across the width of the shell (visible in Figure~\ref{fig:stripOptimised}). 

\subsection{Shell roof over a rectangular domain}
\begin{figure}
  \centering 
  \subfloat[][Initial  geometry ($J (\ary x, \ary u)= 2108$).]
  {\label{fig:roofRefInitial}
  \includegraphics[scale=0.85]{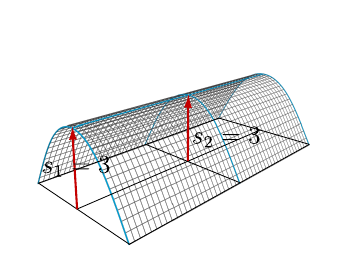}} \\
  \subfloat[][Optimised geometry ($J (\ary x, \ary u)= 149.68$)]
  {\label{fig:roofRefFinal}
  \includegraphics[scale=0.85]{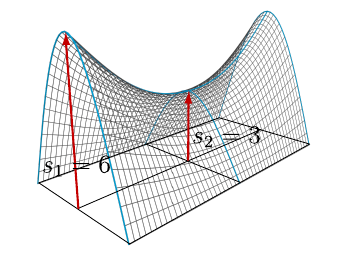}}
\caption{Roof over a rectangular domain.  Initial  and optimised geometries from Bletzinger and Ramm~\cite{bletzinger1993form}. The given compliance cost functions values are computed with the introduced thin-shell solver.}
  \label{fig:roofRef}
\end{figure}

%
Our second example is adapted from Bletzinger et al.~\cite{bletzinger1993form} and considers the compliance optimisation of a roof over a rectangular domain,  see Figure~\ref{fig:roofRef}. The roof has the  plan of~$6 \times 12$ and is pinned along its two long edges.  The applied loading consists of a uniformly distributed vertical load of~$5000$.  The shell thickness is~$t=0.05$; the Young's modulus is~$E=3\times10^{10}$; and the Poisson's ratio is~$\nu = 0.2$. As indicated in Figure~\ref{fig:roofRef},  Bletzinger et al. used only the two height parameters~$s_1$ and~$s_2$ to optimise the roof geometry while maintaining a bi-parabolic shape. Moreover, they chose a cylindrical shell as their initial shape prior to optimisation. 

\begin{figure*}
  \centering 
  {\includegraphics[width=0.9\textwidth]{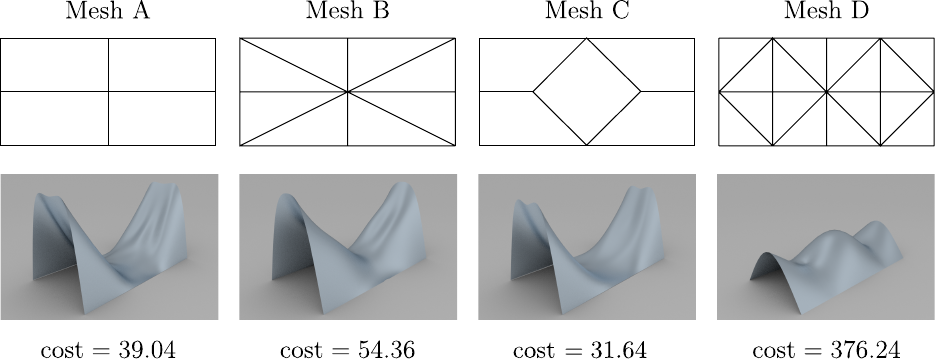}}
  \caption{Roof over a rectangular domain with four different coarse control meshes. Each column shows the coarse control mesh, the obtained optimised geometry and the corresponding cost function value.  The analysis level is chosen with $\ell_c \equiv 2$. The optimisation level starts with $ \ell_o = 0$ and is incremented while $ \ell_o \leq \ell_c$. 
  \label{fig:roofOpt} }
\end{figure*}

In the performed computations a flat rectangle with four different control mesh layouts is used as the initial geometry, see Figure~\ref{fig:roofOpt}. The aim of considering the different initial meshes is to highlight the mesh dependence of the optimised shape. With a gradient-based algorithm a certain mesh dependence of the optimisation results is unavoidable because most structural optimisation problems are non-convex and often not well-posed. In the computed four different cases,  during optimisation only the out-of-plane position of vertices are modified with an prescribed upper bound of~$6$ in order to reproduce the effect of limiting the maximum height~$s_1 \le 6$ in~\cite{bletzinger1993form}. In all cases the analysis level is chosen with~$\ell_c \equiv 2$. The optimisation level starts with~$\ell_o=0$ and is incremented each time the cost function reaches a steady state as long as~$\ell_o \le \ell_c \equiv 2$. In Figure~\ref{fig:roofOpt} the optimised geometries and the corresponding cost function values are shown. The given cost function values are obtained using a fine computational mesh, which has a comparable resolution in all four cases.  For comparison, the cost function of the initial flat rectangular plate is~$\approx 46 \times 10^4$. In all four cases a large reduction in the cost function is achieved and there are significant differences in the final geometries. 
Only starting off with very few optimisation variables in the initial coarse mesh, like in A,  B and C, appears to give lower minima. Note also the resemblance of the optimised geometries for meshes A, B and C to the optimised geometry of Bletzinger et al.~\cite{bletzinger1993form} shown in Figure~\ref{fig:roofRefFinal}. In contrast to the geometry in Figure~\ref{fig:roofRefFinal}, the results for meshes A and B have lower cost function values possibly due to the presence of the fine-scale ripples on the optimised geometries that can be seen in Figure~\ref{fig:roofOpt}, first, second and third columns.

\begin{figure*}
  \centering 
  {\includegraphics[width=0.9\textwidth]{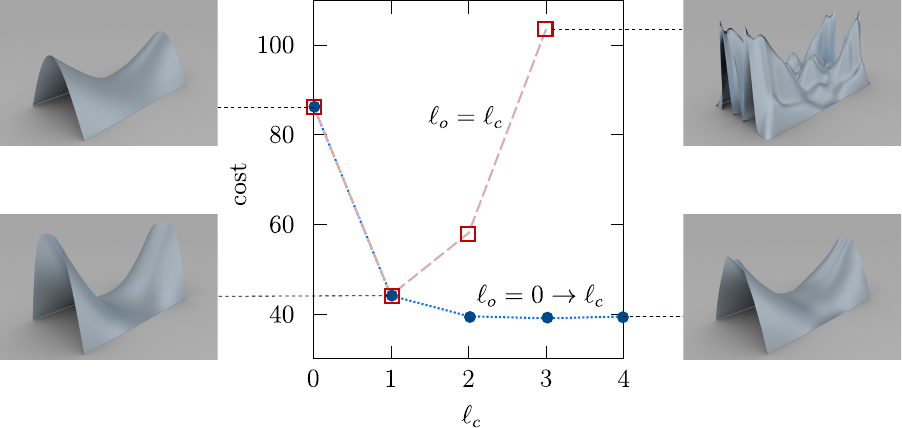}}
  \caption{Comparison of single-resolution and multiresolution optimisation for a roof over a rectangular domain.  The coarse control mesh in all cases is the quadrilateral mesh~A shown in Figure~\ref{fig:roofOpt}.  The dashed red line indicates single resolution optimisation with~$\ell_o \equiv \ell_c$. The dotted blue line indicates multiresolution optimisation, which starts with optimisation level~$\ell_o=0$ and is incremented while~$\ell_o \le \ell_c$. Each data point in the plot denotes the final cost of an independent optimisation study. Inset pictures are the optimised shapes at the indicated data points.   \label{fig:roofMultires} }
\end{figure*}

The choice of the optimisation and computational levels~$\ell_o$ and~$\ell_c$ is  studied next, see Figure~\ref{fig:roofMultires}.  In one set of computations, referred to as {\em single-resolution optimisation}, the two levels are chosen to be the same~$\ell_o \equiv \ell_c$. This means the optimisation variables are simply the out-of-plane positions of the vertices of the computational mesh.  In a second set of computations, referred to as {\em multiresolution optimisation}, the optimisation level always starts with~$\ell_o=0$ and is successively incremented as long as~$\ell_o \le \ell_c$, as discussed in the preceding paragraph. The obtained cost functions for the two different optimisation strategies and for~$\ell_c \in \{0, 1, 2, 3, 4 \}$ are plotted in Figure~\ref{fig:roofMultires}. In the cost function plot each point represents one independent optimisation problem. For all problems the initial control mesh is the quadrilateral mesh~A shown in Figure~\ref{fig:roofOpt}. Note that as expected the result for the multiresolution optimisation for large~$\ell_c \ge 2$ is the same as in Figure~\ref{fig:roofOpt}, first column. In contrast, for single-resolution optimisation the obtained geometries  become highly oscillatory when the level~$\ell_o \equiv \ell_c$ is increased. This is a well-known problem in shape optimisation and suggests the ill-posedness of  the optimisation problem requiring some form of regularisation~\cite{Braibant:1984aa}.  In multiresolution optimisation the design sensitivities are computed on the fine control  mesh and are subsequently projected to the coarser optimisation control mesh using the coarsening matrix~$\ary R$.  As in established  smoothing, or filtering, techniques in shape optimisation~\cite{le2011gradient,bletzinger2014consistent},  the projection of the design sensitivities leads to a smoothing of the design sensitivities. The lack of a projection, and hence a smoothing, in single-resolution optimisation  appears to be the cause of the appearance of the non-optimal jagged geometries with fine-scale oscillations.  In Figure~\ref{fig:roofMultires} all the given cost functions are obtained using a fine control mesh with the same resolution for all data points.

%

\subsection{Freeform architectural shell roof}
%
In design practice, such as in architectural engineering, in addition to structural efficiency there are a number of equally important, often not explicitly quantifiable, design objectives. Although there is a dearth of research on the use of optimisation in a professional design setting,  a recent study shows  that designers usually use optimisation for generating ideas, that is to discover new and unexpected geometries, and do not see  it as a means for generating the ultimate design~\cite{bradner2014parameters}.  With this in mind, isogeometric shape optimisation can  aid the designers in  search for structurally efficient geometries that satisfy all design objectives. An efficient shell structure  can be generated by intermittently shape optimising and manually editing  the control mesh vertex positions and increasing or decreasing the refinement level.  That is, the designer can consult shape optimisation  throughout the entire design stage as often as needed. This sketched design workflow is only feasible with isogeometric analysis and the afforded tight link between the geometry and analysis models.  

\begin{figure*}
	\centering 
	\begin{minipage}[b]{0.45\textwidth}
  		\subfloat[][Perspective view.]{\includegraphics[width=\textwidth]{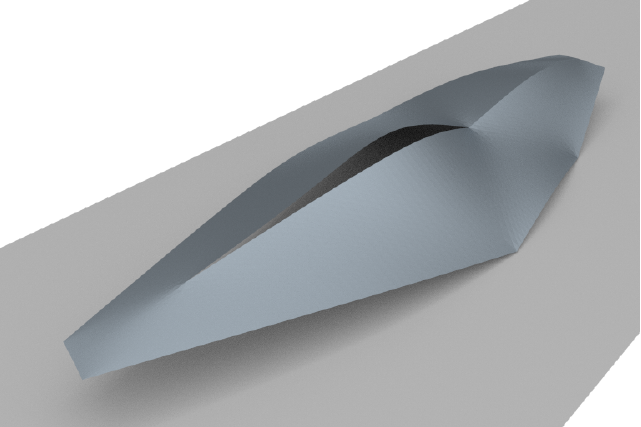}}
	\end{minipage}
		\hspace{0.05\linewidth}
	\begin{minipage}[b][0.29\textwidth][s]{0.375\textwidth}	
		\subfloat[][Front view.]{\includegraphics[width=\textwidth]{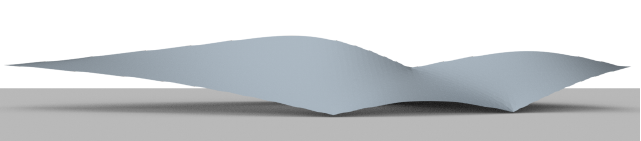}}
		\vfill
		\subfloat[][Back view.]{\includegraphics[width=\textwidth]{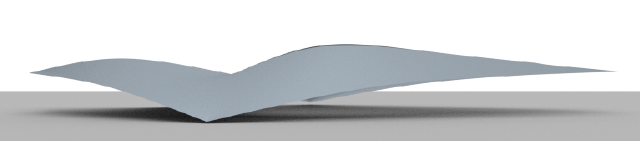}}
	\end{minipage}	
  \caption{Freeform architectural roof supported at three points and containing a central opening and a creased ridge (with $G^0$ continuity). The corresponding coarse resolution control mesh is shown in Figure~\ref{fig:ArchiRoofTags}.}
  \label{fig:ArchiRoofOriginalSubdiv}
\end{figure*}

\begin{figure*}
	\centering 
  	\begin{minipage}[b]{0.45\textwidth}
  		\subfloat[][Top view.]{\includegraphics[width=\textwidth]{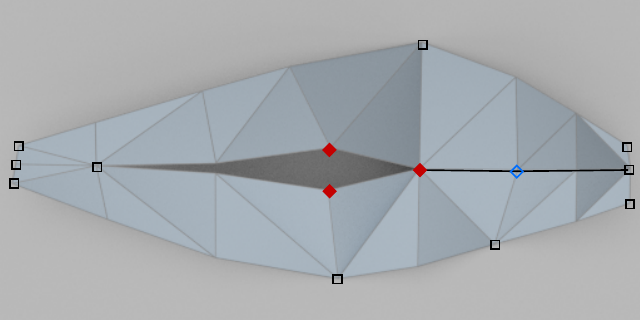}}
	\end{minipage}
	\hspace{0.05\linewidth}
 	\begin{minipage}[b][0.23\textwidth][s]{0.375\textwidth}
		\subfloat[][Front view.]{\includegraphics[width=\textwidth]{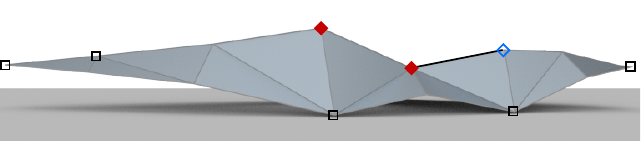}}
		\vfill	
 		\subfloat[][Back view.]{\includegraphics[width=\textwidth]{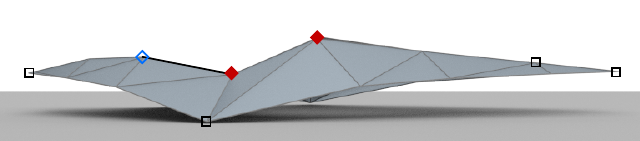}}
	\end{minipage}	
  \caption{Triangular coarse control mesh  with vertex and edge tags. The black square vertices indicate corner tags and the thick black edges indicate crease tags. The effect of the tagging on the limit surface  can be observed in Figure~\ref{fig:ArchiRoofOriginalSubdiv}. During the considered three optimisation design case studies (Design A, Design B and Design C)  some of the vertex coordinates are fixed as specified in Table~\ref{tab:designs}.}
  \label{fig:ArchiRoofTags}
\end{figure*}

To illustrate the use of multiresolution optimisation in a more realistic design setting  we consider the roof structure shown in Figure~\ref{fig:ArchiRoofOriginalSubdiv}. The approximate dimensions of the shell structure are~$2.31 \times 6.27 \times 0.75$.  It is supported only at three points and  contains at the top an opening for lighting purposes and  a crease ($G^0$-continuous feature line).    The applied loading consists of a uniformly distributed load of~$-1000$ (downwards). The shell thickness is~$t=0.02$; the Young's modulus is~$E=10^{10}$; and the Poisson's ratio is~$\nu=0.2$.  The triangular control mesh with~$26$ vertices is depicted in Figure~\ref{fig:ArchiRoofTags}. On the control mesh some of the vertices are tagged as {\em corner vertices} (empty squares) and some of the edges as {\em crease edges} (thick lines).    
  At tagged vertices and edges locally modified stencils  are  applied during subdivision refinement, see~\cite{Biermann:2000aa,Cirak:2011aa} for details. These stencils ensure that sharp corners and $G^0$-continuous feature lines are preserved. The visual effect of the tagging on the limit surface in Figure~\ref{fig:ArchiRoofOriginalSubdiv} is evident. Moreover, the tags have an influence on the entries in the subdivision matrix~$\ary S$, the coarsening matrix~$\ary R$ and the basis functions~$N_i(\theta^1, \theta^2)$.

\begin{figure*}
	\centering 
  	\begin{minipage}[b]{0.45\textwidth}
  		\subfloat[][Perspective view.]{\includegraphics[width=\textwidth]{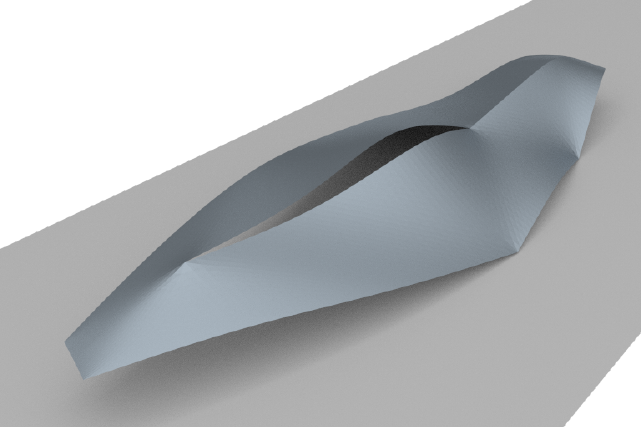}}
	\end{minipage}
	\hspace{0.05\linewidth}
	\begin{minipage}[b][0.29\textwidth][s]{0.375\textwidth}
   		\subfloat[][Front view.]{\includegraphics[width=\textwidth]{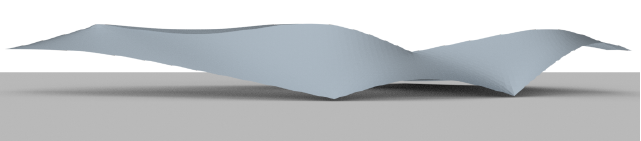}}
  		\vfill
  		\subfloat[][Back view.]{\includegraphics[width=\textwidth]{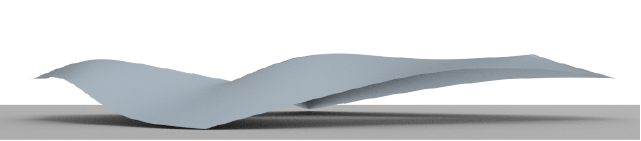}}
	\end{minipage}	
  	\caption{Limit surface of the optimised roof Design A. The final value of the compliance is $23.84$ representing a $23.97\%$ reduction.}
  \label{fig:ArchiRoofsubDesignA}
\end{figure*}

\begin{figure*}
  \centering 
  	\begin{minipage}[b]{0.45\textwidth}
  		\subfloat[][Perspective view.]
  			{\includegraphics[width=\textwidth]{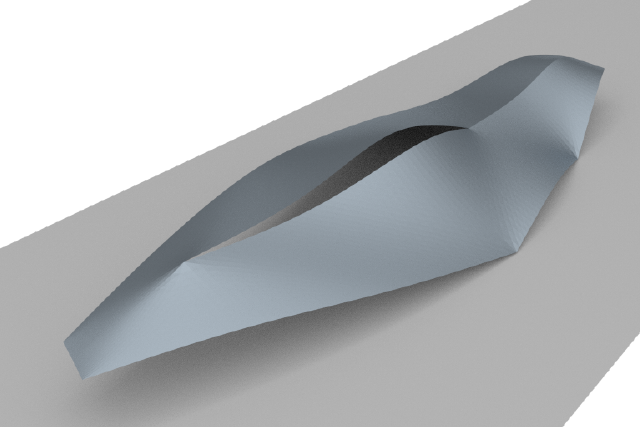}}
  	\end{minipage}
	\hspace{0.05\linewidth}
	\begin{minipage}[b][0.29\textwidth][s]{0.375\textwidth}
  		\subfloat[][Front  view]{\includegraphics[width=\textwidth]{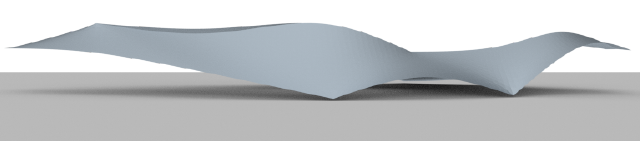}} 
		\vfill
		 \subfloat[][Back view]{\includegraphics[width=\textwidth]{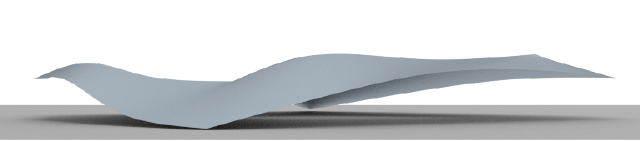}}
	 \end{minipage}
  \caption{Limit surface of the optimised architectural roof  Design B. The final value of the compliance is~$19.17$ representing a~$38.88\%$ reduction.}
  \label{fig:ArchiRoofsubDesignB}
\end{figure*}

\begin{figure*}
	\centering 
  	\begin{minipage}[b]{0.45\textwidth} 
  		\subfloat[][Perspective view.]{\includegraphics[width=\textwidth]{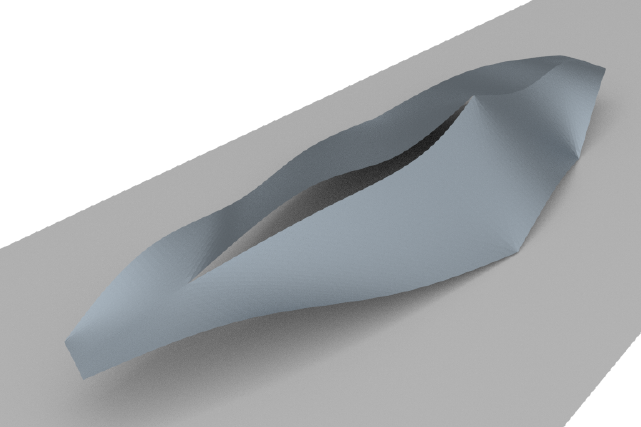}}
	\end{minipage}
	\hspace{0.05\linewidth}
	\begin{minipage}[b][0.29\textwidth][s]{0.375\textwidth}
 		\subfloat[][Front view.]{\includegraphics[width=\textwidth]{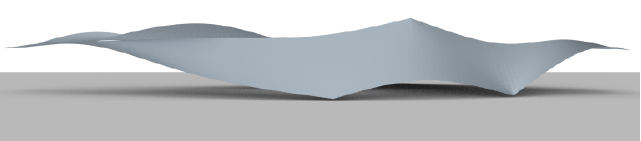}}
 		\hfill 
  		\subfloat[][Back view.]{\includegraphics[width=\textwidth]{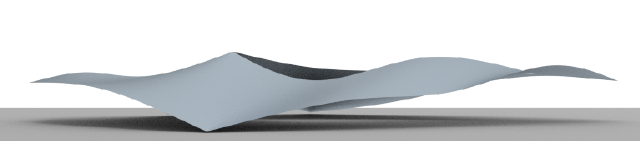}}
	\end{minipage}	
  	\caption{Limit surface of the optimised roof Design C. The final value of the objective function is $6.55$ representing a $79.13\%$ reduction.}
  \label{fig:ArchiRoofsubDesignC}
\end{figure*}

For shape optimisation three different design scenarios, referred to as Design A, Design B and Design C, are considered, see Table~\ref{tab:designs}. Depending on the design scenario the positions of some of the highlighted vertices in Figure~\ref{fig:ArchiRoofTags} are fixed during the optimisation iterations.  
%
\begin{table}
	\centering
	\begin{tabular}{ l |l  } 
		 & Labels of the fixed vertices \\
		 \hline 
		 Design A & square, filled diamond, empty diamond  \\
		 \hline
 		 Design B & square, filled diamond\\ 
		\hline 
		 Design C & square  \\ 
	\end{tabular}
\caption{Three design scenarios for the architectural shell roof. The labels refer to the vertices in Figure~\ref{fig:ArchiRoofTags} whose coordinates are fixed during optimisation. \label{tab:designs}}
\end{table}
 The coarsest control mesh shown in Figure~\ref{fig:ArchiRoofTags} for optimisation contains~$26$ vertices and the twice subdivided finite element mesh with~$\ell_c=2$ contains~$272$ vertices. The design variables in optimisation are  the out-of-plane positions of the vertices. This choice ensures that the planform of the shell roof is preserved.  Only the positions of vertices in levels~$\ell_o=0$ and~$\ell_o=1$ are optimised,  in turn. The optimisation of the vertex positions in the second level results in surfaces with fine-scale oscillations and has been omitted. It is also necessary to restrict the surface area  of the shell to~$A \le 1.2 A_0$, where~$A_0$ is the area of the original surface. The initial value of the compliance cost function is~$31.36$. Design C results in the most efficient optimised shape with a~$79.13\%$ reduction in cost, followed by Design B  with a reduction of~$38.88\%$. The total reduction in the most constrained Design A is only~$23.97\%$. The corresponding optimised shapes for each design scenario are shown in Figures~\ref{fig:ArchiRoofsubDesignA}, \ref{fig:ArchiRoofsubDesignB} and \ref{fig:ArchiRoofsubDesignC}. As can be seen in these figures the optimised Design C has more variation in the surface curvature, which usually leads to stiffer and less compliant structures.  This variation of the curvature is especially pronounced in the large overhanging front part of the shell structure.  Note  also that on all surfaces the characteristic ridge feature and the sharp corners at both ends of the shell are  preserved due to the use of extended subdivision schemes.  
%

\subsection{Conclusions}
%
We introduced the isogeometric shape optimisation of shell structures using triangular Loop and quadrilateral Catmull-Clark subdivision surfaces as a common representation for geometric modelling and finite element analysis. The presented examples demonstrate that  efficient and flexible representation of  freeform shell geometries is essential for shape optimisation. In the implemented gradient-based shape optimisation approach more optimal geometries are obtained when, starting from the coarsest control mesh, the vertices of increasingly finer control meshes are chosen as geometric design variables.  Irrespective of the control mesh resolution for optimisation a sufficiently fine control mesh can always be used for finite element discretisation. In general the finite element control meshes have to be finer than the optimisation control meshes  for accuracy reasons.  The introduced approach effectively allows the designer to choose an optimal geometry with a visually pleasing and technically feasible smallest feature size. With the increasing availability of subdivision surfaces in CAD systems it is expected that it will become feasible to import the optimised geometries back into a CAD environment for continuing with the design process.

The presented isogeometric shape optimisation approach can be extended and improved in several ways. First,  we considered only the structural compliance as a  cost function and that for only one loading case. In practice, there are many more load cases and  competing cost functions which have to be taken into account. For instance, the structural stability, i.e. buckling,  of  optimised thin shells is often critical and has to be taken into account~\cite{reitinger1995buckling}.  Moreover, the geometric fidelity of the surfaces in the presented approach can be improved using more  recent higher-degree subdivision surfaces, such as  the NURBS-compatible subdivision surfaces~\cite{Cashman:2009aa}, or adaptively refined subdivision surfaces~\cite{wei2015truncated, Bornemann:2013aa}. Finally, the multiresolution editing techniques from computer graphics based on wavelet-like decomposition of surfaces can be employed for interlacing geometry generation by the user with automated structural optimisation~\cite{bandara2016shape}.


\bibliographystyle{elsarticle-num-names}
\bibliography{shellOptimisation}

\end{document}